\documentclass[preprint,review,10pt]{elsarticle}
\usepackage{amsfonts}
\usepackage{amsmath,amssymb}
\usepackage{graphicx}
\usepackage{setspace}
\usepackage[left=2cm,top=2cm,right=2cm]{geometry}

\newtheorem{lemma}{Lemma}

\newtheorem{remark}{Remark}

\newtheorem{theorem}{Theorem}
\numberwithin{equation}{section} \journal{XYZ}
\begin{document}
\title{Rate of convergence of certain families of Jain operators of integral type}
\author[label*,label1,label2]{Prashantkumar Patel}
\ead{prashant225@gmail.com}
\author[label1,label3]{Vishnu Narayan Mishra}
\ead{vishnu\_narayanmishra@yahoo.co.in; vishnunarayanmishra@gmail.com}
\address[label1]{Department of Applied Mathematics \& Humanities,
S. V. National Institute of Technology, Surat-395 007 (Gujarat), India}
\address[label2]{Department of Mathematics,
St. Xavier's College(Autonomous), Ahmedabad-380 009 (Gujarat), India}
\address[label3]{L. 1627 Awadh Puri Colony Beniganj, Phase-III, Opposite - Industrial Training Institute (ITI), Ayodhya Main Road, Faizabad, Uttar Pradesh
224 001, India}
\fntext[label*]{Corresponding authors}
\begin{abstract}
In the present paper, the authors introduce and investigate new
sequences of positive linear operators which include some well
known operators as special cases. Here we estimate the rate of
convergence for functions having derivatives of bounded variation
by families of Jain operators of integral type.
\end{abstract}
\begin{keyword} Positive linear operators; Jain operators; Bounded variation; Rate of convergence\\
\textit{2000 Mathematics Subject Classification: } primary 41A25, 41A30, 41A36. \end{keyword}

\maketitle
\section{Introduction}
In the year 1972, Jain \cite{jain1972approximation} introduced
and studied the new class of positive linear operators using Poisson-type distribution as
\begin{eqnarray}\label{RPS2.ch0.2}
G_n^{\mu}(f,x) = \sum_{k=0}^{\infty} n x(nx+k\mu)^{k-1}\frac{e^{-(nx+k\mu)}}{k!}f\left(\frac{k}{n}\right),
\end{eqnarray}
where $\mu\in [0,1)$ and $f\in C(\mathbf{R}^{+})$. In the particular case $\mu=0$, $G_n^{0}$, $n\in \mathbf{N}$,
 turn out to be well-known Sz\'{a}sz-Mirakjan operators \cite{szasz1950generalization}.
 Umar and Razi \cite{umar1985approximation} studied Kantorovich-type extensions of $G_n^{\mu}$.
 Tarabie \cite{tarabie2012jain} and
 Mishra and Patel \cite{mishrasome2013}  introduced integral versions of the Jain operators using Beta
 basis functions and discussed their approximation properties. Recently, both the authors have
 established the Jain-Baskakov operators and different generalizations of them in \cite{patel2015jain}.\\
\indent The general integral modification of Jain-Baskakov
operators to approximate Lebesgue integrable functions on the
interval $[0,\infty)$, can be defined as follows:
\begin{equation}\label{11.1.eq1.1}
K_{n}^{\mu,r}(f,x)=\frac{n^r(n-r-1)!}{(n-2)!} \sum_{v=0}^{\infty} \omega_{\mu}(v,n x)\int_0^{\infty} p_{n-r,v+r}(t)f(t) dt, ~~r\geq 0,
\end{equation}
where $n\in \mathbb{N}$, $ r\in \mathbb{N} \cup\{0\}$, $n>r$,
$\mu\in[0,1)$ and the Baskakov and the Jain basis functions are
defined as
$$ p_{n,v}(t) = {n+v-1 \choose v}\frac{t^{v}}{(1+t)^{n+v}},$$
$$\omega_{\mu}(v,nx) = nx(nx+v\mu)^{v-1}\frac{e^{-(nx+v\mu)}}{v!}.$$
We note that, the operators $K_{n}^{\mu,r}(f,\cdot)$ are linear
and positive. If $r=0$, the operators \eqref{11.1.eq1.1} are equal
to the Jain-Baskakov operators studied by Patel and Mishra
\cite{patel2015jain}. The rate of convergence for functions having
derivatives of bounded variation is the investigated by many
authors
\cite{gupta1999rate,srivastava2003certain,zeng2002rate,zeng1998rate,patel2014rate,zeng2000rate}.
Here we extend these studies to investigation of similar
properties of the modified Jain-Baskakov operators as defined in
\eqref{11.1.eq1.1}.
\section{Estimation of moments}
\begin{lemma}[\cite{jain1972approximation}]\label{13.1.lemma1}
For $G_{n}^{\mu} (t^m,x),~~m=0,1,2$, we have
$$ G_{n}^{\mu} (1,x)= 1,~~~ G_{n}^{\mu} (t,x)= \frac{x}{1-\mu},~~~ G_{n}^{\mu} (t^2,x)= \frac{x^2}{(1-\mu)^2}+ \frac{x}{n(1-\mu)^3}.$$
\end{lemma}
\begin{lemma}\label{12.1.lemma3}
Let the $\textrm{m}^{\textrm{th}}$ order moment be defined as
\begin{equation}
K_{n,r,m}^{\mu}(x) = (n-r-1) \sum_{v=0}^{\infty} \omega_{\mu}(v,n x)\int_0^{\infty} p_{n-r,v+r}(t)t^m dt.
\end{equation}
Then $$ K_{n,r,0}^{\mu}(x) = 1, n>r+1, ~~
K_{n,r,1}^{\mu}(x)=\frac{n x + (r+1)
(1-\mu)}{(n-r-2)(1-\mu)}\text{ for } n>r+2,$$
 $$K_{n,r,2}^{\mu}(x)=\frac{1}{(n-r-2)(n-r-3)}\left[\frac{n^2  x^2}{(1-\mu)^2}+ \left[\frac{n}{(1-\mu)^3}
            + \frac{n(2r+3)}{1-\mu}\right]x + (r+1)(r+2)\right] \text { for } n>r+3.$$
\end{lemma}
\textbf{Proof: } For $n>r+1$, we have
\begin{eqnarray*}
K_{n,r,0}^{\mu}(x) &=& (n-r-1) \sum_{v=0}^{\infty} \omega_{\mu}(v,n x)\frac{(n+v-1)!}{(v+r)!(n-r-1)!}\int_0^{\infty} \frac{t^{v+r}}{(1+t)^{n+v}} dt\\
&=& (n-r-1) \sum_{v=0}^{\infty} \omega_{\mu}(v,n x)\frac{(n+v-1)!}{(v+r)!(n-r-1)!}\cdot \frac{(n-r-2)!(v+r)!}{(n+v-1)!} = G_{n}^{\mu} (1,x)= 1.
\end{eqnarray*}
Now, for $m=1$ and $n>r+2$, we get
\begin{eqnarray*}
K_{n,r,1}^{\mu}(x) &=& (n-r-1) \sum_{v=0}^{\infty} \omega_{\mu}(v,n x)\frac{(n+v-1)!}{(v+r)!(n-r-1)!}\int_0^{\infty} \frac{t^{v+r}}{(1+t)^{n+v}}\cdot t  dt\\
&=& (n-r-1) \sum_{v=0}^{\infty} \omega_{\mu}(v,n x)\frac{(n+v-1)!}{(v+r)!(n-r-1)!}\cdot \frac{(n-r-3)!(v+r+1)!}{(n+v-1)!} \\
&=& (n-r-1) \sum_{v=0}^{\infty} \omega_{\mu}(v,n x)\frac{v+r+1}{(n-r-1)(n-r-2)}\\
&=& \frac{n}{(n-r-2)} \sum_{v=0}^{\infty} \omega_{\mu}(v,n x) \frac{v}{n}+ \frac{r+1}{(n-r-2)}\sum_{v=0}^{\infty} \omega_{\mu}(v,n x)\\
&=& \frac{n}{(n-r-2)} \frac{x}{(1-\mu)}+ \frac{r+1}{(n-r-2)}\\
&=& \frac{n x + (r+1) (1-\mu)}{(n-r-2)(1-\mu)}.
\end{eqnarray*}
Further, for $m=2$ and $n> r+3$, we obtain
\begin{eqnarray*}
K_{n,r,2}^{\mu}(x)&=& (n-r-1) \sum_{v=0}^{\infty} \omega_{\mu}(v,n x)\frac{(n+v-1)!}{(v+r)!(n-r-1)!}\int_0^{\infty} \frac{t^{v+r}}{(1+t)^{n+v}}\cdot t^2  dt\\
&=& (n-r-1) \sum_{v=0}^{\infty} \omega_{\mu}(v,n x)\frac{(n+v-1)!}{(v+r)!(n-r-1)!}\cdot \frac{(n-r-4)!(v+r+2)!}{(n+v-1)!} \\
&=&  (n-r-1) \sum_{v=0}^{\infty} \omega_{\mu}(v,n x)\frac{(v+r+1)(v+r+2)}{(n-r-1)(n-r-2)(n-r-3)}\\
&=& \frac{1}{(n-r-2)(n-r-3)} \sum_{v=0}^{\infty} \omega_{\mu}(v,n x)(v^2 + (2r+3)v+(r+1)(r+2))\\
&=& \frac{1}{(n-r-2)(n-r-3)} \left[n^2 G_{n}^{\mu} (t^2,x) + n(2r+3) G_{n}^{\mu} (t,x) + (r+1)(r+2) G_{n}^{\mu} (1,x)\right]\\
&=& \frac{1}{(n-r-2)(n-r-3)}\left[\frac{n^2  x^2}{(1-\mu)^2}+
\left[\frac{n}{(1-\mu)^3}+ \frac{n(2r+3)}{1-\mu}\right]x +
(r+1)(r+2)\right].
\end{eqnarray*}
This completes the proof of Lemma \ref{12.1.lemma3}.
\begin{lemma}\label{11.1.lemma1}
Let the $\textrm{m}^{\textrm{th}}$ order central moment be defined
as
\begin{equation}
T_{n,r,m}^{\mu}(x) = (n-r-1) \sum_{v=0}^{\infty} \omega_{\mu}(v,n x)\int_0^{\infty} p_{n-r,v+r}(t)(t-x)^m dt
\end{equation}
then 
$T_{n,r,0}^{\mu}(x) =1$ for $n>r+1$, $\displaystyle
T_{n,r,1}^{\mu}(x) = \frac{(1+r)(1-\mu) +x ((2+r)(1- \mu) +n \mu
)}{(n-r-2)(1-\mu)}\textrm{ for } n>r+2 $ and
\begin{eqnarray*}
T_{n,r,2}^{\mu}(x) &=& x^2 \left(1+\frac{n^2}{(n-r-3) (n-r-2)
(1-\mu )^2}-\frac{2 n}{(n-r-2)
(1-\mu )}\right)\\
&&+x \left(\frac{n \left(1+(3+2 r) (1-\mu )^2\right)}{(n-r-3)
(n-r-2) (1-\mu )^3}-\frac{2 (1+r)}{n-r-2}\right)\\
&&+\frac{(1+r) (2+r)}{(n-r-3) (n-r-2)}, \textrm{ for } n>r+3.
\end{eqnarray*}
\end{lemma}
The proof of above lemma follows from linear properties of the
operators $K_{n}^{\mu,r}(f,\cdot)$. Lemma \ref{12.1.lemma3} shows
that the operators $K_{n}^{\mu,r}(f,\cdot)$ does not preserve the
linear functions, that is, $K_{n}^{\mu,r}(f,x)\neq f(x) $ for
$f(t)= at +b$, where $a$ and $b$ are real constants.

\begin{remark}\label{12.1.remark1}
For fixed $r$, there is a constant $C > 1$ (which depends only on
$r$) such that, for all $n$ being sufficiently large, all $\mu$
being sufficiently small (say $0 <\mu<\mu_0$), and $x \in
(0,\infty)$,
$$\frac{x^2}{nC}\leq T_{n,r,2}^{\mu}(x) \leq \frac{Cx^2}{n}.$$
\end{remark}
\begin{remark}
By using Cauchy-Schwarz inequality, it follows from Remark
\ref{12.1.remark1}, that for all $n$ being sufficiently large, all
$\mu$ being sufficiently small, $C > 1$ and $x\in (0,\infty)$, we
have
\begin{equation}
(n-r-1) \sum_{v=0}^{\infty} \omega_{\mu}(v,n x)\int_0^{\infty}
p_{n-r,v+r}(t)|t-x| dt\leq
\left[T_{n,r,2}^{\mu}(x)\right]^{1/2}\leq \sqrt{\frac{Cx^2}{n}}.
\end{equation}
\end{remark}
\begin{lemma}\label{12.1.lemma2}
Let $x\in (0,\infty)$ and $C>1$, then for all $n$ being
sufficiently large and all $\mu$ being sufficiently small, we have
\begin{equation}
\delta_{n,r}^{\mu}(x,y) = (n-r-1) \sum_{v=0}^{\infty}
\omega_{\mu}(v,n x)\int_0^{y} p_{n-r,v+r}(t) dt\leq
\frac{Cx^2}{n(x-y)^2},~~0\leq x < y
\end{equation}
\begin{equation}
1- \delta_{n,r}^{\mu}(x,z) = (n-r-1) \sum_{v=0}^{\infty}
\omega_{\mu}(v,n x)\int_z^{\infty} p_{n-r,v+r}(t) dt\leq
\frac{Cx^2}{n(z-x)^2}, ~~x<z<\infty.
\end{equation}
\end{lemma}
\textbf{Proof: } The proof of the above lemma follows easily by
using Remark \ref{12.1.remark1}. For instance, for the first
inequality for all $n$ being sufficiently large, all $\mu$ being
sufficiently small and $0 \leq  y < x$, we have
\begin{eqnarray*}
\delta_{n,r}^{\mu}(x,y) &=& (n-r-1) \sum_{v=0}^{\infty} \omega_{\mu}(v,n x)\int_0^{y} p_{n-r,v+r}(t) dt\\
&\leq & (n-r-1) \sum_{v=0}^{\infty} \omega_{\mu}(v,n x)\int_0^{y} p_{n-r,v+r}(t) \frac{(t-x)^2}{(y-x)^2}dt\\
&=& \frac{T_{n,r,2}^{\mu}(x)} {(y-x)^2}\leq \frac{Cx^2}{n(x-y)^2}.
\end{eqnarray*}
The proof of the second inequality follows along the similar
lines.

\section{Rate of Convergence}
By $DB_q (0,\infty)$ (where $q$ is some positive integer), we
mean the class of absolutely continuous functions $f$ defined on
$(0,\infty)$ satisfying the following conditions:
\begin{itemize}
\item[(i)] $f (t) = O(t^q ),~ t \to \infty$;
\item[(ii)] the function $f$ has the first derivative on the interval $(0,\infty)$ which coincide
almost everywhere with a function which is of bounded variation on
every finite subinterval of $(0,\infty)$. It can be observed that
for all functions $f \in DB_q (0,\infty)$, we can have the
representation
$$f(x) = f(c) + \int_c^x \psi(x) dt, ~~~~0< c \leq x.$$
\end{itemize}
\begin{theorem}\label{12.1.theorem1}
Let $f \in  DB_q (0,\infty)$, $q > 0$ and $x \in (0,\infty)$.
Then for $C > 1$, all $\mu$  being sufficiently small and all $n$
being sufficiently large, we have
\begin{eqnarray*}
\bigg|  \frac{(n-2)!}{n^r (n-r-2)!}K_{n}^{\mu,r}(f,x)-f(x) \bigg|
& \leq &\frac{C x}{n}\left(\sum_{k=1}^{[\sqrt{n}]}
                        \bigvee_{x-x/k}^{x+x/k} (( f')_x) + \frac{x}{\sqrt{n}}\bigvee_{x-x/\sqrt{n}}^{x+x/\sqrt{n}} (( f')_x)\right)\\
&& +\frac{C}{n}\left(| f (2x) - f (x) - x f'(x^{+})| + | f (x)|\right)\\
&& + |f'(x^{+})|\frac{C x }{n} + \frac{1}{2}\sqrt{\frac{C x^2}{n}} | f'(x^+) - f'(x^-) |\\
&&+ \frac{1}{2} | f'(x^+)+ f'(x^-)| \frac{(1+r)(1-\mu) +x ((2+r)(1- \mu) +n \mu )}{(n-r-2)(1-\mu)}+ O(n^{-q}),
\end{eqnarray*}
where $\bigvee_a^b f (x)$ denotes the total variation of $f_x$ on $[a, b]$, and $f_x$ is defined by
\begin{equation*}f_x(t) = \left\{
  \begin{array}{lc}
    f(t)-f(x^-),  & 0\leq t < x; \\
     0, & t=x; \\
    f(t) - f(x^{+}), & x<t < \infty. \\
  \end{array}
\right.
\end{equation*}
\end{theorem}
\textbf{Proof: }
Using the mean value theorem, we can write
\begin{eqnarray*}
\bigg|  \frac{(n-2)!}{n^r (n-r-2)!}K_{n}^{\mu,r}(f,x)-f(x) \bigg|&
\leq &
(n-r-1) \sum_{v=0}^{\infty} \omega_{\mu}(v,n x)\int_0^{\infty} p_{n-r,v+r}(t)|f(t) -f(x) | dt\\
&=& \int_0^{\infty} \bigg| \int_x^t (n-r-1) \sum_{v=0}^{\infty} \omega_{\mu}(v,n x)p_{n-r,v+r}(t)f'(u)du\bigg| dt.
\end{eqnarray*}
Also, using the identity
\begin{eqnarray*}
f'(u) &=& \frac{f'(x^+) +f'(x^-)}{2} + (f')_x(u) +
\frac{f'(x^+)-f'(x^-)}{2}~\textrm{sgn}(u-x) + \left[f'(x)
-\frac{f'(x^+) +f'(x^-)}{2} \right]\chi_x(u),
\end{eqnarray*}
where
\begin{equation*}
\chi_x(u)=\left\{
\begin{array}{lc}
1, & u=x;\\
0, & u\neq x.\end{array}\right.
\end{equation*}
Obviously, we have
\begin{eqnarray*}
(n-r-1) \int_0^{\infty} \left(\int_x^t \left[f'(x) -\frac{f'(x^+)
+f'(x^-)}{2} \right]\chi_x(u) du \right)
                            \sum_{v=0}^{\infty} \omega_{\mu}(v,n x)p_{n-r,v+r}(t)dt=0.
\end{eqnarray*}
Thus, using above identities, we can write
\begin{eqnarray}\label{12.1.eq3.1}
&&\bigg|  \frac{(n-2)!}{n^r (n-r-2)!}K_{n}^{\mu,r}(f,x)-f(x) \bigg|\nonumber\\
& \leq & \int_0^{\infty} \bigg|
                \int_x^t (n-r-1) \sum_{v=0}^{\infty} \omega_{\mu}(v,n x)p_{n-r,v+r}(t)\left(\frac{f'(x^+)+f'(x^-)}{2} + (f')_x(u)\right)du\bigg| dt\nonumber\\
&&+ \int_0^{\infty}
            \bigg| \int_x^t (n-r-1) \sum_{v=0}^{\infty} \omega_{\mu}(v,n x)p_{n-r,v+r}(t)\left(\frac{f'(x^+)-f'(x^-)}{2}~\textrm{sgn}(u-x)\right)du\bigg| dt.
\end{eqnarray}
Also, it can be verified that
\begin{eqnarray}\label{12.1.eq3.2}
&& \bigg| (n-r-1) \int_0^{\infty} \left(\int _x^t \frac{f'(x^+)-f'(x^-)}{2}~\textrm{sgn}(u-x) du \right)
                        \sum_{v=0}^{\infty} \omega_{\mu}(v,n x)p_{n-r,v+r}(t) dt\bigg|\nonumber\\
&\leq&\bigg| \frac{f'(x^+)-f'(x^-)}{2}\bigg|
\left[T_{n,r,2}^{\mu}(x)\right]^{1/2}
\end{eqnarray}
and
\begin{eqnarray}\label{12.1.eq3.3}
 \bigg| (n-r-1) \int_0^{\infty} \left(\int _x^t \frac{f'(x^+)+ f'(x^-)}{2} du \right)
                        \sum_{v=0}^{\infty} \omega_{\mu}(v,n x)p_{n-r,v+r}(t)dt\bigg|= \bigg| \frac{f'(x^+)+f'(x^-)}{2}\bigg| T_{n,r,1}^{\mu}(x).
\end{eqnarray}
Combining \eqref{12.1.eq3.1}-\eqref{12.1.eq3.3}, we have
\begin{eqnarray}\label{12.1.eq3.4}
\bigg|  \frac{(n-2)!}{n^r (n-r-2)!}K_{n}^{\mu,r}(f,x)-f(x) \bigg|
    &\leq& \bigg| (n-r-1) \int_x^{\infty} \left(\int _x^t (f')_x(u)du \right)\sum_{v=0}^{\infty} \omega_{\mu}(v,n x)p_{n-r,v+r}(t)dt\bigg|\nonumber\\
&&+\bigg| (n-r-1) \int_0^{x} \left(\int _x^t (f')_x(u)du \right)\sum_{v=0}^{\infty} \omega_{\mu}(v,n x)p_{n-r,v+r}(t)dt\bigg|\nonumber\\
&&+ \bigg| \frac{f'(x^+)-f'(x^-)}{2}\bigg| \left[T_{n,r,2}^{\mu}(x)\right]^{1/2} + \bigg| \frac{f'(x^+)+f'(x^-)}{2}\bigg| T_{n,r,1}^{\mu}(x)\nonumber\\
&=& |A_{n,r}^{\mu}(f,x)| + |B_{n,r}^{\mu}(f, x)| \nonumber\\
&&+  \bigg| \frac{f'(x^+)-f'(x^-)}{2}\bigg|
\left[T_{n,r,2}^{\mu}(x)\right]^{1/2} + \bigg|
\frac{f'(x^+)+f'(x^-)}{2}\bigg| T_{n,r,1}^{\mu}(x).
\end{eqnarray}
Applying Remark \ref{12.1.remark1} and Lemma \ref{11.1.lemma1}, in \eqref{12.1.eq3.4}, we have
\begin{eqnarray}\label{12.1.eq}
\bigg|  \frac{(n-2)!}{n^r (n-r-2)!}K_{n}^{\mu,r}(f,x)-f(x)\bigg|
                    &\leq& |A_{n,r}^{\mu}(f,x)| + |B_{n,r}^{\mu}(f,x)|
                        +  \bigg| \frac{f'(x^+)-f'(x^-)}{2}\bigg|\sqrt{\frac{Cx^2}{n}}  \nonumber\\
&&+ \bigg| \frac{f'(x^+)+f'(x^-)}{2}\bigg|\frac{(1+r)(1-\mu) +x
((2+r)(1- \mu) +n \mu )}{(n-r-2)(1-\mu)}.
\end{eqnarray}
In order to complete the proof of the theorem it sufficient to
estimate the terms $A_{n,r}^{\mu}(f,x)$ and $B_{n,r}^{\mu}(f,x)$.
            Applying integration by parts and Lemma \ref{12.1.lemma2}
with $\displaystyle y= x- \frac{x}{\sqrt{n}}$, we have
\begin{eqnarray}
| B_{n,r}^{\mu}(f,r) | &=& \bigg|(n-r-1) \int_0^{x} \left(\int _x^t (f')_x(u)du d_t (\delta_{n,r}^{\mu} (x,t) ) \right)\bigg|\nonumber\\
&=& \bigg| \int_0^{x}  \delta_{n,r}^{\mu} (x,t) (f')_x(t) dt\bigg|\nonumber\\
&\leq & \left( \int_0^y + \int_y^x \right) |(f')_x(t)|| \delta_{n,r}^{\mu} (x,t)| dt\nonumber\\
&\leq &  \frac{C x^2}{n} \int_0^y \bigvee_t^x ((f')_x)\frac{1}{(t-x)^2}dt + \int_y^x \bigvee_t^x\left( (f')_x\right)dt \nonumber\\
&\leq &  \frac{C x^2}{n} \int_0^y \bigvee_t^x\left(
(f')_x\right)\frac{1}{(t-x)^2}dt
                + \frac{x}{\sqrt{n}} \bigvee_{x-\frac{x}{\sqrt{n}}}^x \left((f')_x\right).\nonumber
\end{eqnarray}
Let $\displaystyle u= \frac{x}{x-t}$. Then we have
\begin{eqnarray*}
\frac{C x^2}{n} \int_0^y \bigvee_t^x\left(
(f')_x\right)\frac{1}{(t-x)^2}dt
        &=& \frac{C x^2}{n} \int_1^{\sqrt{n}} \bigvee_{x-\frac{x}{u}}^x \left((f')_x\right) du\\
        &\leq& \frac{C x}{n} \sum_{k=1}^{[\sqrt{n}]} \bigvee_{x-\frac{x}{k}}^x \left((f')_x\right).
\end{eqnarray*}
Thus
\begin{eqnarray}\label{12.1.eq3.5}
| B_{n,r}^{\mu}(f,r) | &\leq& \frac{C x}{n}\sum_{k=1}^{[\sqrt{n}]}
                                        \bigvee_{x-\frac{x}{k}}^x\left(
                                        (f')_x\right)
                                        + \frac{x}{\sqrt{n}} \bigvee_{x-\frac{x}{\sqrt{n}}}^x \left((f')_x\right).
\end{eqnarray}
On the other hand, we have
\begin{eqnarray}\label{12.1.eq3.6}
|A_{n,r}^{\mu}(f,x)| &=& \bigg| (n-r-1) \int_x^{\infty} \left(\int _x^t (f')_x(u)du \right)\sum_{v=0}^{\infty} \omega_{\mu}(v,n x)p_{n-r,v+r}(t)dt\bigg|\nonumber\\
&=& \bigg| (n-r-1) \int_{2x}^{\infty} \left(\int _x^t (f')_x(u)du
\right)\sum_{v=0}^{\infty} \omega_{\mu}(v,n
x)p_{n-r,v+r}(t)dt\nonumber\\
            &&+ \int_x^{2x} \left(\int _x^t (f')_x(u)du \right)d_t(1 -\delta_{n,r}^{\mu}(x,t) )\bigg|\nonumber\\
&\leq & \bigg| (n-r-1)  \sum_{v=0}^{\infty} \omega_{\mu}(v,n x) \int_{2x}^{\infty} (f(t) -f(x)) p_{n-r,v+r}(t)dt\bigg|\nonumber\\
&&+|f'(x^+)|\bigg| (n-r-1)  \sum_{v=0}^{\infty} \omega_{\mu}(v,n x) \int_{2x}^{\infty} (t-x) p_{n-r,v+r}(t)dt\bigg|\nonumber\\
&& + \bigg| \int_x^{2x} (f')_x(u) du \bigg|| 1- \delta_{n,r}^{\mu}(x,2x)| + \int_x^{2x} |(f')_x(t)|| 1- \delta_{n,r}^{\mu}(x,t)|dt \nonumber\\
&\leq &(n-r-1)  \sum_{v=0}^{\infty}\omega_{\mu}(v,n x) \int_{2x}^{\infty} Mt^{2q} p_{n-r,v+r}(t)dt\nonumber\\
&& +\frac{|f(x)|}{x^2}(n-r-1)  \sum_{v=0}^{\infty}\omega_{\mu}(v,n x) \int_{2x}^{\infty} p_{n-r,v+r}(t)(t-x)^2 dt\nonumber\\
&& +|f'(x^+)|(n-r-1)  \sum_{v=0}^{\infty}\omega_{\mu}(v,n x) \int_{2x}^{\infty} p_{n-r,v+r}(t)|t-x| dt\nonumber\\
&& + \frac{C}{n} | f (2x) - f (x) - x f'(x^+)|\nonumber\\
&&+ \frac{C x}{n} \sum_{k=1}^{[\sqrt{n}]}
\bigvee_{x}^{x+\frac{x}{k}}\left((f')_x\right)+
\frac{x}{\sqrt{n}}\bigvee_{x}^{x+\frac{x}{\sqrt{n}}}\left((f')_x\right).
\end{eqnarray}
To estimate the integral $\displaystyle (n-r-1)
\sum_{v=0}^{\infty}\omega_{\mu}(v,n x) \int_{2x}^{\infty} Mt^{2q}
p_{n-r,v+r}(t)dt$, in \eqref{12.1.eq3.6} above, we proceed as
follows: Obviously $t \geq 2x$  implies that $t \leq  2(t- x)$ and
it follows from Lemma \ref{11.1.lemma1}, that
\begin{eqnarray*}
(n-r-1)  \sum_{v=0}^{\infty}\omega_{\mu}(v,n x) \int_{2x}^{\infty} Mt^{2q} p_{n-r,v+r}(t)dt
&\leq& M 2^{2q} (n-r-1)  \sum_{v=0}^{\infty}\omega_{\mu}(v,n x) \int_{0}^{\infty}  p_{n-r,v+r}(t)(t-x)^{2q}dt\\
&=& M2^{2q}T_{n,r,2q}^{\mu}(x)  = O(n^{-q}). (n \to \infty)
\end{eqnarray*}
Applying Schwarz inequality and Remark \ref{12.1.remark1}, third term in right hand side of \eqref{12.1.eq3.6} is estimated as follows:
\begin{eqnarray*}
&&|f'(x^+)|(n-r-1)  \sum_{v=0}^{\infty}\omega_{\mu}(v,n x)
\int_{2x}^{\infty} p_{n-r,v+r}(t)|t-x| dt\\
& \leq &\frac{
|f'(x^+)|}{x} (n-r-1)
            \sum_{v=0}^{\infty}\omega_{\mu}(v,n x) \int_{0}^{\infty} p_{n-r,v+r}(t)(t-x)^2 dt\\
&=& \frac{C x}{n}|f'(x^+)|.
\end{eqnarray*}
Thus by Lemma \ref{11.1.lemma1} and Remark \ref{12.1.remark1}, we have
\begin{eqnarray}\label{12.1.eq3.8}
|A_{n,r}^{\mu} ( f, x)| &\leq& O(n^{-q}) + \frac{C x}{n} |f'(x^+)|
+ \frac{C }{n}
                         \left(| f (2x) - f (x) - x f'(x^+)|+ |f(x)|\right)\nonumber\\
&&+ \frac{C x}{n} \sum_{k=1}^{[\sqrt{n}]}
\bigvee_{x}^{x+\frac{x}{k}}((f')_x)
                        + \frac{x}{\sqrt{n}}\bigvee_{x}^{x+\frac{x}{\sqrt{n}}}((f')_x).
\end{eqnarray}
Collecting the estimates \eqref{12.1.eq}, \eqref{12.1.eq3.5} and \eqref{12.1.eq3.8}, we get the required result.\\
This completes the proof of Theorem \ref{12.1.theorem1}.
\section{Modification of the operators $K_{n}^{\mu,r}$ with parameter $c$}
In the year 1995, Gupta \textit{et al.} \cite{gupta1995onsim}
introduced integral modification of the Sz\'{a}sz-Mirakyan
operators by considering the weight functions of Beta basis
functions. Recently, Dubey and Jain \cite{dubey2008rate} modified
the operators discussed by Gupta   \textit{et al.}
\cite{gupta1995onsim}  with a parameter $c > 0$ and studied their
approximation properties. This type of approach was also discussed
by many authors, we refer some papers as
\cite{gupta2009appro,gupta2011note,mishrasome2013,verma2012convergence}.
This motivated us to study the rate of convergence for the
generalized Jain-Baskakov operators \eqref{11.1.eq1.1} with
parameter $c$, which is defined as the follows:
\begin{equation}\label{11.1.eq1.2}
K_{n,c}^{\mu,r}(f,x)=\frac{n^r\Gamma(\frac{n}{c}-r)}{\Gamma(\frac{n}{c}-1)} \sum_{v=0}^{\infty}
                 \omega_{\mu}(v,n x)\int_0^{\infty} p_{n-rc,v+r}(t,c)f(t) dt, ~~r\geq 0,
\end{equation}
where $n\in \mathbb{N}$, $n>rc$, the generalized Baskakov basis
function defined as $\displaystyle  p_{n,v}(t,c) =
\frac{\Gamma(\frac{n}{c}
+v)}{\Gamma(\frac{n}{c})\Gamma(v+1)}\frac{(ct)^{v}}{(1+ct)^{\frac{n}{c}+v}}$
    and $\omega_{\mu}(v,nx)$ as defined in \eqref{11.1.eq1.1}. If $c=1$,
            then the operators \eqref{11.1.eq1.2} reduce to  the operators defined in \eqref{11.1.eq1.1}.
\begin{lemma}\label{13.2.lemma1.2}
Let the $\textrm{m}^{\textrm{th}}$ order moment be defined as
\begin{equation}
K_{n,r,m}^{\mu,c}(x) = (n-rc-c) \sum_{v=0}^{\infty} \omega_{\mu}(v,n x)\int_0^{\infty} p_{n-rc,v+r}(t)t^m dt.
\end{equation}
Then $$ K_{n,r,0}^{\mu,c}(x) = 1, ~~  K_{n,r,1}^{\mu,c}(x)=\frac{n
x + (r+1)(1-\mu)}{(n- c r-2c)(1-\mu)},\textrm{ for } n>(r+2)c,$$
$$K_{n,r,2}^{\mu,c}(x)=\frac{1}{(n-c r-2c)(n-c
r-3c)}\left[\frac{n^2  x^2}{(1-\mu)^2}+ \left[\frac{n}{(1-\mu)^3}+
\frac{n(2r+3)}{1-\mu}\right]x + (r+1)(r+2)\right], \text{ for }
n>(r+3)c.$$
\end{lemma}
The proof of above Lemma follows along the lines of Lemma
\ref{12.1.lemma3}; thus, we omit the details.
\begin{lemma}\label{11.2.lemma1}
Let the $\textrm{m}^{\textrm{th}}$ order moment be defined as
\begin{equation}
T_{n,r,m}^{\mu}(x,c) = (n-rc-c) \sum_{v=0}^{\infty} \omega_{\mu}(v,n x)\int_0^{\infty} p_{n-rc,v+r}(t,c)(t-x)^m dt
\end{equation}
then $T_{n,r,0}^{\mu}(x,c) =1$ , $\displaystyle T_{n,r,1}^{\mu}(x,c)
                    =  \frac{(1+r)(1-\mu) +x (c (2 + r)(1-\mu) +n \mu  )}{(n- c r-2c)(1-\mu)},\text { for } n>(r+2)c$ and
\begin{eqnarray*} T_{n,r,2}^{\mu}(x,c)
    &=& x^2 \left(1+\frac{n^2}{(n-2c-r c) (n-3c-r c) (1-\mu )^2}-\frac{2 n}{(n-2c -r c) (1-\mu )}\right)\\
    &&+x \left(\frac{n \left(1+(3+2 r) (1-\mu )^2\right)}{(n-2c-r c) (n-3c- r c) (1-\mu )^3}-\frac{2 (1+r)}{n-2c-r c}\right)\\
    &&+\frac{(1+r) (2+r)}{(n-2c-r c) (n-3c-r c)} , \text { for }
n>(r+3)c.
\end{eqnarray*}
\end{lemma}
\begin{remark}\label{12.2.remark1}
For fixed $r$, there is a constant $C_1> 1$ (which depends only on
$r$) such that, for all $n$ being sufficiently large, all $\mu$
being sufficiently small (say $0 <\mu<\mu_0$), and $x \in
(0,\infty)$,
$$\frac{x^2}{C_1n}\leq T_{n,r,2}^{\mu}(x,c) \leq \frac{C_1x^2}{n}.$$
\end{remark}
Let $B_{x^2}[0,\infty)$ = \{ $f$ : for every $x\in [0,\infty), |f(x)|\leq M_f(1+x^2), M_f$ being
a constant depending on $f$\}. By $C_{x^2}[0,\infty)$, we denote the subspace of all continuous
functions belonging to $B_{x^2}[0,\infty)$. Also, $C_{x^2}^{*}[0,\infty)$ is  subspace of
all function $f\in C_{x^2}[0,\infty)$ for which $\displaystyle \lim_{x\to \infty} \frac{f(x)}{1+x^2}$ is finite.
The norm on $C_{x^2}^{*}[0,\infty)$ is $\displaystyle \|f\|_{x^2} = \sup_{x\in [0,\infty)} \frac{|f(x)|}{1+x^2}$.\\
Examining relation given in Lemma \ref{13.2.lemma1.2} and based on famous Korovkin theorem \cite{korovkin1953convergence},
 it is clear that $\{K_{n,c}^{\mu,r}\}, n>rc$ does not form an approximation process. To enjoying of this property,
 we replace the constant $\mu$ by a number $\mu_n\in [0, 1)$ with $$\displaystyle \lim_{n\to \infty}\mu_n = 0.$$
Now by Lemma 2.1, we ensure that $$\lim_{n\to \infty}
K_{n,c}^{\mu_n,r} (t^m, x) = x^m ,~~~m = 0, 1, 2$$ uniformly on
any interval compact $E\subset[0,\infty)$. Base on Korovkin's
criteria,  we can state the following result:
\begin{theorem} Let $K_{n,c}^{\mu_n,r}(f,x)$ with $n > rc > 0$, be defined as in \eqref{11.1.eq1.2},
where $\displaystyle \lim_{n\to \infty} \mu_n = 0$. For any
compact set $E \subset [0,\infty)$ and for each $f \in
C_{x^2}^{*}[0,\infty)$ one has
$$\displaystyle \lim_{n \to \infty}K_{n,c}^{\mu_n,r}(f,x) = f(x), \textrm{ uniformly in }  x \in E.$$
\end{theorem}
\begin{theorem}
Let $f \in  DB_q (0,\infty)$, $q > 0$ and $x \in (0,\infty)$.
Then for $C_1>1$, all $\mu$ being sufficiently small and $n$ being
sufficiently large, we have
\begin{eqnarray*}
\bigg|  \frac{\Gamma(\frac{n}{c}-1)}{n^r \Gamma(\frac{n}{c}-r-1)}\left(K_{n,c}^{\mu,r}(f,x)-f(x)\right) \bigg|
    & \leq &\frac{C_1x}{n}\left(\sum_{k=1}^{[\sqrt{n}]} \bigvee_{x-x/k}^{x+x/k} (( f')_x)
            + \frac{x}{\sqrt{n}}\bigvee_{x-x/\sqrt{n}}^{x+x/\sqrt{n}} (( f')_x)\right)\\
&& +\frac{C_1}{n}\left(| f (2x) - f (x) - x f'(x^{+})| + | f (x)|\right)\\
&& + |f'(x^{+})|\frac{C_1x }{n} + \frac{1}{2}\sqrt{\frac{C_1x^2}{n}} | f'(x^+) - f'(x^-) |+ O(n^{-q})\\
&&+ \frac{1}{2} | f'(x^+)+ f'(x^-)| \frac{(1+r)(1-\mu) +x (c (2 +
r)(1-\mu) +n \mu  )}{(n- c r-2c)(1-\mu)},
\end{eqnarray*}
where $\bigvee_a^b f (x)$ denotes the total variation of $f_x$ on
$[a, b]$, and $f_x$ is defined by
\begin{equation*}f_x(t) = \left\{
  \begin{array}{lc}
    f(t)-f(x^-),  & 0\leq t < x; \\
     0, & t=x; \\
    f(t) - f(x^{+}), & x<t < \infty. \\
  \end{array}
\right.
\end{equation*}
\end{theorem}
The proof of the above theorem follows along the lines of Theorem \ref{12.1.theorem1}; thus we omit the details.
\begin{remark}
In \cite{acar2011rate} Acar \textit{et al.}, estimated the rate of
convergence for functions having derivatives of bounded variation
in simultaneous approximation of  the general integral
modification of the Sz\'{a}sz -Mirakyan operators
 having the weight functions of Baskakov basis functions. It would be interesting to extend study of the operators
 \eqref{11.1.eq1.1} in simultaneous approximation. To achieve simultaneous approximation,
 one has to  establish recurrence relation in terms of derivative of the Jain operators $G_{n}^{\mu}(f,x)$.
 We are still curious to know, what is recurrence relation of moments of the Jain operators?
\end{remark}
\begin{remark}
In year 1983, Stancu \cite{Stancu1983Appr} generalized the
Bernstein polynomials with two parameters $\alpha$ and $\beta$
($0\leq \alpha\leq \beta)$. After this many classical
approximation operators generalized into two parameters $\alpha$
and $\beta$
\cite{verma2012some,mishra2013approximation2,mishra2014durrmeyer,atakut2010stancu,acar2014simultaneous,mishra2013some,mishra2013inverse,mishra2013statistical}.\\
A Stancu generalization of the operators \eqref{11.1.eq1.1}, is
introduced as follows:
\begin{equation}\label{11.1.eq1.3}
K_{n,\alpha,\beta}^{\mu,r}(f,x)=\frac{n^r(n-r-1)!}{(n-2)!} \sum_{v=0}^{\infty}
             \omega_{\mu}(v,n x)\int_0^{\infty} p_{n-r,v+r}(t)f\left(\frac{nt+\alpha}{n+\beta}\right) dt, ~~r\geq 0,
\end{equation}
where $\omega_{\mu}(v,n x)$ and $p_{n,v}(t)$ as defined in
\eqref{11.1.eq1.1}. For the operators  \eqref{11.1.eq1.3}, one can
study its local approximation properties,
 Voronovskaja type asymptotic results and rate of convergence for functions having derivatives of bounded variation.
\end{remark}
\section*{Acknowledgements}
The authors would like express their deep gratitude to the
anonymous learned referee(s) which resulted in the subsequent
improvement of this research article.

\end{document}